\documentclass[12pt]{article}

\usepackage{amsmath,amsthm,amssymb}
\allowdisplaybreaks
\usepackage{mathrsfs}
\usepackage{cite}

\usepackage[pdftex,bookmarksopen=true,bookmarks=true,unicode,setpagesize]{hyperref}
\hypersetup{colorlinks=true,linkcolor=black,citecolor=black}

\newtheorem{theorem}{Theorem}

\theoremstyle{remark}

\theoremstyle{remark}

\theoremstyle{remark}
\newtheorem{remark}[theorem]{Remark}

\newcommand{\R}{\mathbb R}
\newcommand{\di}{\partial}
\begin{document}

\vspace{-20mm}
\begin{center}{\large \bf
The projection spectral theorem and Jacobi fields\\[4mm]
{\it Dedicated to Yuri Makarovych Berezansky\\
on the occasion of his 90th birthday}}
\end{center}

{\large Eugene Lytvynov}\\ Department of Mathematics,
Swansea University, Singleton Park, Swansea SA2 8PP, U.K.;
e-mail: \texttt{e.lytvynov@swansea.ac.uk}

{\small
\begin{center}
{\bf Abstract}
\end{center}
\noindent We review several applications of Berezansky's projection spectral theorem to Jacobi fields in a symmetric Fock space, which lead to L\'evy white noise measures.} \vspace{2mm}

Running head: Jacobi fields

AMS (2010) Subject Classification: 60H40, 60G57, 47B36

\section{Introduction}

The projection spectral theorem for a family of commutating self-adjoint (or more generally, normal) operators and the corresponding Fourier transform in generalized joint eigenvectors of this family play a fundamental role in functional analysis and its applications, in particular, to infinite dimensional (stochastic) analysis and mathematical physics. In the case where the operators family has additionally a cyclic vector, these result allow  to realize the operators from the family as multiplication operators   in an $L^2$-space with respect to a probability measure. While this theory has been studied by quite a few authors, including Yu.~Berezansky, L. G\r{a}rding, I. Gelfand, A.~Kostyuchechko, G. Kats, and K. Maurin, in the most general setting, this result was proved by Berezansky in \cite{ber_1984}, see also \cite{ber_1985,ber_1978} and Chapter 3 in \cite{BK}.
In particular, the theorems of Berezansky allow the operators family to be uncountable, which is a very natural assumption for many applications.

In the simplest case where there is just one self-adjoint operator, this theory is deeply connected with the theory of orthogonal polynomials on the real line, see Chapter VII in \cite{ber_1968}. More precisely,  one starts with a Hermitian operator $J$ in $\ell_2$ defined on all finite vectors. The operator $J$ is assummed to be defined by an infinite Jacobi (i.e., tridiagonal) matrix. Under proper conditions on this matrix, the operator $J$ is essentially self-adjoint and its closure, $\widetilde J$, can be realized (through a unitary isomorphism) as an operator of multiplication in a space $L^2(\R,d\mu)$, where $\mu$ is a probability measure on $\R$.
Furthermore, the set of polynomials  is dense in $L^2(\R,d\mu)$, while the elements of the Jacobi matrix $J$ are precisely the coefficients of the recursive formula satisfied by the orthonormal polynomials in $L^2(\R,d\mu)$.

This theory admits important extensions to an infinite dimensional setting. More precisely, instead of $\ell_2$, one considers a  Hilbet space $\mathcal F$  of the form
\begin{equation}\label{hugf7r}
\mathcal F=\bigoplus_{n=0}^\infty\mathcal F_n,\end{equation}
and a family of (generally speaking unbounded) commuting self-adjoint operators $A(\varphi)$, $\varphi\in\Phi$, where  $\Phi$ is an index set (in many applications, $\Phi$ is a nuclear space). The operators $A(\varphi)$ are assumed to have a tridiagonal structure with respect to the orthogonal decomposition \eqref{hugf7r}. Such a family of operators, $(A(\varphi))_{\varphi\in\Phi}$, is called a (commutative) Jacobi field, see \cite{BKosh}.
An application of the projection spectral theorem to the family
$(A(\varphi))_{\varphi\in\Phi}$ leads to a unitary isomorphism between $\mathcal F$ and an $L^2$ space with respect to a probability measure $\mu$ (if $\Phi$ is a nuclear space, then typically $\mu$ is a probability measure on $\Phi'$, the dual space of $\Phi$). This approach was first proposed by Berezansky in \cite{ber_1993}, see also \cite{BLL}.

In this paper, we will review three examples of application of the projection spectral theorem to  commutative Jacobi fields, which lead to important classes of probability measures on an infinite dimensional space.
 All these measures will be L\'evy white noise measures, or more generally, generalized stochastic processes with independent values in the sense of \cite{GV}.

The paper is organized  as follows. In Section~\ref{cte64ue}, we will formulate the projection spectral theorem in the form which is convenient for our studies. In Section~\ref{ghfdy7er5}, we will show how the  projection spectral theorem  leads to the standard Gaussian measure. In Section~\ref{huf7tr}, we will discuss the case of Poisson measure, and in Section~\ref{gyd6ei6} we will discuss the case of a L\'evy white noise measure. (The latter case will generalize the results of Sections~\ref{ghfdy7er5} and \ref{huf7tr}.) Finally, in Section~\ref{vgfvyf} we will shortly discuss further developments related to the application of the projection spectral theorem and free (noncommutative) L\'evy white noise .

\section{The projection spectral theorem}\label{cte64ue}

Let us first recall the spectral theorem in the case of one self-adjoint operator. Let $H$ be a real, separable Hilbert space and let  $(A, D(A))$ be a self-adjoint operator.
Let $\Omega \in H $ and assume that $ \Omega $ is cyclic for $A$, i.e., $\Omega \in D(A^n)$, $n \in \mathbb{N}$, and the linear span of the set $\{ \Omega, A\Omega, A^2\Omega, \dots \}$ is  dense in $H$. Then, the spectral theorem implies that there exists a unique probability measure $\mu$ on $\mathbb{R}$ such that the  mapping  $I$ given by  $$(I\Omega)(x) = 1, \quad
(IA^n\Omega)(x) = x^n, \quad x\in\R,\ n \in \mathbb{N},$$ extends by linearity and continuity to a unitary operator $$ I: H \rightarrow L^2(\mathbb{R}, \mu).$$
Furthermore, $IAI^{-1}$ is  the operator
of multiplication by $x$ in $L^2(\R,\mu)$.
In fact, the measure $\mu$ is given by
$$\mu(\alpha) = (E(\alpha)\Omega, \Omega)_H, \quad \alpha \in \mathcal{B}(\mathbb{R}),$$ where $E(\cdot)$ is the resolution
of the identity of $A$ and $\mathcal{B}(\mathbb{R})$ is the Borel $\sigma$-algebra on $\R$. The measure $\mu$ is called the spectral measure of $A$  (at $\Omega$).

The following theorem generalizes the above result to the case of a family of commuting self-adjoint operators indexed by elements of a nuclear space.
This theorem  follows immediately from Chapter 3 of\cite{BK}.

\begin{theorem} \label{spectral}
Assume that we have two standard triples
$$\Phi \subset H \subset \Phi'\quad\text{and}\quad \Psi \subset \mathcal{F} \subset \Psi' ,$$ where $H$ and $\mathcal{F}$
are real separable Hilbert spaces, $\Phi$ and $\Psi$ are nuclear spaces, which are densely embedded into $H$ and $\mathcal F$, respectively, and $\Phi'$ ($\Psi'$ respectively) is the dual space of $\Phi$ ($\Psi$ respectively) with respect to the center space $H$ ($\mathcal F$ respectively). For $\omega\in\Phi'$ and $\varphi\in\Phi$, we denote by $\langle\omega,\varphi\rangle$ the dual pairing between $\omega$ and $\varphi$.

Assume that we have a family $(A(\varphi))_{\varphi \in \Phi} $ of Hermitian operators in $\mathcal{F}$ which satisfy the following conditions:

\begin{enumerate}
\item $D(A(\varphi )) = \Psi$, $\varphi \in\Phi $.

\item $A(\varphi )\Psi \subset \Psi$ for each $\varphi \in \Phi$,
and furthermore $ A(\varphi ):\Psi \rightarrow \Psi $ is continuous.

\item $A(\varphi_1)A(\varphi_2)f = A(\varphi_2)A(\varphi_1)f$,  $f \in \Psi$ (i.e., the operators $A(\varphi )$ algebraically commute on $\Psi $).

\item For all $f, g \in \Psi$, the mapping  $$ \Phi \ni \varphi \mapsto (A(\varphi )f,g)_{\mathcal{F}} \in \mathbb{R}$$  is continuous.

\item There exists a vector $\Omega$ in $\mathcal{F}$ which is cyclic for $(A(\varphi ))_{\varphi \in \Phi }$, i.e., the linear span of the    set
$$\{\Omega,\, A(\varphi_1) \cdots A(\varphi_k)\Omega \mid \varphi_1, \dots, \varphi_k \in \Phi , \ k \in \mathbb{N}\}$$ is dense in $\mathcal {F}$.

\item for any $f \in \Psi $ and  $\varphi \in \Phi $, the vector $f$ is analytic for the operator $A(\varphi )$, i.e., for some $t>0$,
$$\sum_{n=0}^\infty t^n\,\frac{\|A(\varphi)f\|_{\mathcal F}}{n!}<\infty.$$
\end{enumerate}

Then, each operator $A(\varphi )$, $\varphi \in\Phi $, is essentially self-adjoint on $\Psi$ and we denote its closure by $(\tilde{A}(\varphi ),D(\tilde{A}(\varphi )))$.
These operators commute in the sense of their resilutions of the identity. Furthermore, there exists a unique probability measure
$\mu$ on $(\Phi',\mathcal{C}(\Phi'))$ such that the linear operator $I: \mathcal{F} \rightarrow L^2(\Phi',\mu)$ given by
$I \Omega = 1$ and
\begin{align*}
I(\tilde{A}(\varphi_1) \cdots \tilde{A}(\varphi_n)\Omega) &= I(A(\varphi_1) \cdots A(\varphi_n)\Omega)\\&= \langle  \cdot, \varphi_1 \rangle \cdots \langle \cdot ,\varphi_n \rangle
\in L^2(\Phi', \mu)
\end{align*}
is unitary. Here, $\mathcal C(\Phi')$ denotes the cylinder $\sigma$-algebra on $\Phi'$.

Under the action of $I$, each  operator $(\tilde{A}(\varphi),D(\tilde{A}(\varphi)))$, $\varphi\in\Phi$, becomes the operator of multiplication by
$\langle \cdot ,\varphi \rangle$ in $L^2(\Phi', \mu)$, denoted by $M(\varphi)$ and given by
$$ D(M(\varphi))=\{F \in L^2(\Phi', \mu): \langle\cdot, \varphi\rangle F \in L^2(\Phi', \mu)\}$$
and for each $F\in  D(M(\varphi))$,
$$ (M(\varphi) F)(\omega)=\langle\omega,\varphi\rangle F(\omega).$$
\end{theorem}

\section{Gaussian measure}\label{ghfdy7er5}

Let us now discuss how a  standard Gaussian measure appears as a result of application of the projection spectral theorem, see \cite{BLL} for details and proofs.

For a real separable Hilbert space $H$, we denote by $\mathcal F(H)$ the (real) symmetric Fock space over $H$:
$$\mathcal F(H)=\bigoplus_{n=0}^\infty H^{\odot n}n!\,.$$
Here $\odot$ denotes symmetric tensor product, and for a Hilbert space $\mathfrak H$ and a positive constant $\alpha$, $\mathfrak H\alpha$ denotes the Hilbert space which coincides with $\mathfrak H$ as a set, and with scalar product
$(g,h)_{\mathfrak H\alpha}=(g,h)_{\mathfrak H}\,\alpha$. We also used the notation $H^{\odot 0}=\R$.

Let $\Phi \subset H \subset \Phi'$ be a standard triple as in Theorem~\ref{spectral}. Since $\Phi$ is a nuclear space, it admits a representation
$$\Phi=\operatornamewithlimits{proj\,lim}_{\tau\in T}H_\tau,$$
where $(H_\tau)_{\tau\in T}$ is a family of Hilbert spaces which is directed by embedding (i.e., for any $\tau,\tau_2\in T$ there exists $\tau_3\in T$ such that  $H_{\tau_2}$ is densely and continuously embedded both into $H_{\tau_1}$ and into $H_{\tau_2}$), and for any $\tau_1\in T$ there exists $\tau_2\in T$ such that $H_{\tau_2}\subset H_{\tau_1}$ and the operator of embedding of $H_{\tau_2}$ into $H_{\tau_1}$ is of Hilbert--Schmidt type. Note that
$$\Phi'=\operatornamewithlimits{ind\,lim}_{\tau\in T}H_{-\tau},$$
where $H_{-\tau}$ is the dual space of $H_{\tau}$ with respect to the center space $H$.

Then
$$\Phi^{\odot n}=\operatornamewithlimits{proj\,lim}_{\tau\in T}H_\tau^{\odot n}$$
is called the $n$-th symmetric tensor power of $\Phi$, $\Phi^{\odot n}$ being a nuclear space. We thus get the standard  triple
$$\Phi^{\odot n}\subset H^{\odot n}\subset\Phi'^{\,\odot n},$$
where
$$\Phi'^{\,\odot n}=\operatornamewithlimits{ind\,lim}_{\tau\in T}H_{-\tau}^{\odot n}.$$

We denote by $\mathcal F_{\text{fin}}(\Phi)$ the topological direct sum of the $\Phi^{\odot n}$ spaces. $\mathcal F_{\text{fin}}(\Phi)$ is a nuclear space which consists of all sequences $(f^{(n)})_{n=0}^\infty$ with $f^{(n)}\in\Phi^{\odot n}$ and for some $N\in\mathbb N$, we have  $f^{(n)}=0$, $n\ge N$. Convergence in this space means uniform finiteness of non-zero components and coordinate-wise convergence. Thus, we can construct the standard triple
\begin{equation}\label{huyf67e6}
\mathcal F_{\text{fin}}(\Phi)\subset \mathcal F(H)\subset \mathcal F^*_{\text{fin}}(\Phi), \end{equation}
where $\mathcal F^*_{\text{fin}}(\Phi)$ is the dual space of
 $\mathcal F_{\text{fin}}(\Phi)$ with respect to the center  space $\mathcal F(H)$. The space $\mathcal F^*_{\text{fin}}(\Phi)$ consists of all sequences $(F^{(n)})_{n=0}^\infty$ with $F^{(n)}\in \Phi'^{\,\odot n}$ and convergence in this space means coordinate-wise convergence in the respective $\Phi'^{\,\odot n}$ spaces.

For each $\varphi\in\Phi$ we define a creation operator $a^+(\varphi)$ and an annihilation operator $a^-(\varphi)$ on  $\mathcal F_{\text{fin}}(\Phi)$ by
\begin{align*}
a^+(\varphi) f^{(n)}&=\varphi\odot f^{(n)},\quad \varphi\in\Phi^{\odot n},\\
a^-(\varphi)h^{\otimes n}&=n(\varphi,h)_H\, h^{\otimes (n-1)},\quad h\in\Phi.
\end{align*}
The following theorem was proved in \cite{BLL}.

 \begin{theorem}\label{buf7o6} For each $\varphi\in\Phi$, we set
 $$A(\varphi)=a^+(\varphi)+a^-(\varphi).$$
 Then conditions of Theorem~\ref{spectral} are satisfied with $\mathcal F=\mathcal F(H)$ and $\Psi=\mathcal F_{\text{fin}}(\Phi)$. The corresponding measure $\mu$ is the standard Gaussian measure on $\Phi'$:
 $$\int_{\Phi'}e^{i\langle\omega,\varphi\rangle }\,d\mu(\omega)=\exp\left[-\frac12\|\varphi\|_H^2\right],\quad \varphi\in\Phi.$$

  \end{theorem}

 Let us consider the special case of a Gaussian white noise. Let $X$ be a smooth Riemannian manifold and let $dx$ be the volume measure on it. Let $H=L^2(X,dx)$ and let $\Phi=\mathcal D(X)$ be the space of all smooth functions on $X$ with compact support. $\mathcal D(X)$ can be represented as a projective limit of weighted Sobolev spaces on $X$, which shows, in particular, that $\mathcal D(X)$ is a nuclear space.
 For $x\in X$, we denote by $\di_x$ and $\di_x^\dag$ the  annihilation and creation  operators at point $x$, respectively. Thus,
 $$\di_xf^{(n)}=nf^{(n)}(x,\cdot),\quad \di_x^\dag f^{(n)}=\delta_x\odot f^{(n)},\quad f^{(n)}\in\mathcal D(X)^{\odot n},$$ where $\delta_x\in\mathcal D'(X)$ is the delta function at $x$. We then have
\begin{equation}\label{tyr7i5} A(\varphi)=\int_X \varphi(x)(\di^\dag_x+\di_x)\,dx,\quad \varphi\in\mathcal D(X).\end{equation}
 The measure $\mu$ from Theorem~\ref{buf7o6} is then called a Gaussian white noise measure, and the operators $\di^\dag_x+\di_x$ ($x\in X$) can be thought of as a Gaussian white noise (see e.g.\ \cite{HKPS} for details).

 \section{Poisson measure}\label{huf7tr}
 Let us again assume that $H=L^2(X,dx)$ and $\Phi=\mathcal D(X)$. For $\varphi\in\mathcal D(X)$, we define a neutral operator $a^0(\varphi)$ as the differential second quantization of the operator of multiplication by $\varphi$. Thus, $a^0(\varphi)$ is the linear  operator on $\mathcal F_{\text{fin}}(\mathcal D(X))$ defined by
 $$ a^0(\varphi)f^{(n)}(x_1,\dots,x_n)=\big(\varphi(x_1)+\dots+\varphi(x_n)\big)f^{(n)}(x_1,\dots,x_n),\quad f^{(n)}\in\mathcal D(X)^{\odot n}.$$
 Let now
 \begin{equation}\label{vytd65r}
  A(\varphi)=a^+(\varphi)+a^0(\varphi)+a^-(\varphi)+\int_X\varphi(x)\,dx.\end{equation}
 It is easy to see that
 $$ a^0(\varphi)=\int_X\varphi(x)\di^\dag_x\di_x\,dx. $$
 Hence, analogously to \eqref{tyr7i5}, we get
 \begin{equation}\label{vgufy7} A(\varphi)=\int_X \varphi(x)(\di^\dag_x+\di^\dag_x\di_x+\di_x+1)\,dx,\quad \varphi\in\mathcal D(X).\end{equation}

 The following theorem was proved in \cite{L1}, see also \cite{ber_2000}.

 \begin{theorem}\label{gfye6}
 For each $\varphi\in\mathcal D(X)$, let $A(\varphi)$ be defined by \eqref{vytd65r}. Then conditions of Theorem~\ref{spectral} are satisfied with $\mathcal F=\mathcal F(H)$ and $\Psi=\mathcal F_{\text{fin}}(\Phi)$. The corresponding measure $\mu$ is the Poisson measure on $\mathcal D'(X)$:
 $$\int_{\mathcal D'(X)}e^{i\langle\omega,\varphi\rangle }\,d\mu(\omega)=\exp\left[
 \int_X(e^{i\varphi(x)}-1)\,dx
 \right],\quad \varphi\in\mathcal D(X).$$
 \end{theorem}

According to \eqref{vgufy7} and Theorem~\ref{gfye6}, the operators
$$\di^\dag_x+\di^\dag_x\di_x+\di_x+1=(\di^\dag_x+1)(\di_x+1),\quad x\in X,$$ can be thought of as Poisson white noise, see e.g.\ \cite{IK} for details.

Let $\Gamma(X)$ denote the configuration space over $X$:
$$\Gamma(X)=\{\gamma\subset X\mid \text{for each compact $\Lambda\subset X$, $\gamma\cap \Lambda$ is a finite set }\}. $$
One usually identifies each $\gamma\in\Gamma(X)$ with the Radon measure $\sum_{x\in\gamma}\delta_x$. (Here $\delta_x$ denotes the Dirac measure with mass at $x$.) This gives the inclusion $\Gamma(X)\subset \mathcal D'(X)$. It can be shown that $\Gamma(X)$ is a measurable set in $\mathcal D'(X)$, and furthermore the trace $\sigma$-algebra of $\mathcal C(\mathcal D'(X))$ on $\Gamma(X)$ coincides with the Borel $\sigma$-algebra on $\Gamma(X)$ generated by the vague topology on $\Gamma(X)$. It can be shown that $\mu(\Gamma(X))=1$, so that $\mu$ is  (the distribution of) a Poisson point process in $X$ with intensity measure $dx$, e.g.\  \cite{Kallenberg}.

\begin{remark}More generally, let us fix a parameter $\lambda>0$ and consider the operators
$$(\di^\dag_x+\sqrt\lambda)(\di_x+\sqrt\lambda)=\sqrt\lambda\di^\dag_x+\di^\dag_x\di_x+\sqrt\lambda\di_x+\lambda,\quad x\in X,$$
and the corresponding smeared operators $(A(\varphi))_{\varphi\in\mathcal D(X)}$. The application of Theorem~\ref{spectral}  to these operators yields a Poisson point process in $X$ with intensity measure $\lambda\,dx$:
$$\int_{\Gamma(X)}e^{i\langle\omega,\varphi\rangle }\,d\mu(\omega)=\exp\left[
 \int_X(e^{i\varphi(x)}-1)\lambda\,dx
 \right],\quad \varphi\in\mathcal D(X).$$
\end{remark}

\section{L\'evy white noise}\label{gyd6ei6}

Assume that  $\sigma$ is a finite (nonzero) measure on
$(\mathbb{R},\mathcal{B}(\mathbb{R}))$.
We also assume that the measure $\sigma$ satisfies the estimate
\begin{align} \label{equ 3.2}
\int_\mathbb{R} |s|^n \,d\sigma(s) \leq C^n n!, \quad  n \in
\mathbb{N},
\end{align}
for some $C>0$. Equivalently, there exists an $\varepsilon>0$ such that $\int_\R e^{\varepsilon|s|}\,d\sigma(s)<\infty$.

Let $\mathcal P$ denote the set of polynomials on $\R$.
One can easily introduce a nuclear space topology on $\mathcal P$ which yields the following convergence: $p_n\to p$  in $\mathcal P$ as $n\to\infty$ if and only if the degree of all polynomials $p_n$ is bounded by some finite constant and, for each $k\in\mathbb N_0$, we have $a_{nk}\to a_k$ as $n\to\infty$, where $p_{n}(s)=\sum_{k\ge0}a_{nk}s^k$ and $p(s)=\sum_{k\ge0} a_ks^k$.

If the measure $\sigma$ has an infinite number of points in its support, we can consider $\mathcal P$ as a subset of $L^2(\R,\sigma)$, and due to estimate \eqref{equ 3.2} $\mathcal P$ is a dense subset of $L^2(\R,\sigma)$. Furthermore, $\mathcal P$ is continuously embedded into $L^2(\R,\sigma)$. If the support of the measure $\sigma$ is finite, we may instead use the factorization of $\mathcal P$ with respect to the set $$\{p\in\mathcal P\mid p=0\text{ $\sigma$-a.e.}\}.$$
We will still use $\mathcal P$ to denote this factorization. Note that in this case, $\mathcal P$ is just the set of all real-valued functions on the support of measure $\sigma$.

 Thus, in any case, we get a nuclear space $\mathcal P$ which is topologically  (i.e., densely and continuously) embedded into $L^2(\R,\sigma)$.
We then consider the nuclear space
$$\mathcal D=\mathcal D'(X)\otimes \mathcal P,$$
where $\mathcal D'(X)$ is the nuclear space as in Section~\ref{ghfdy7er5}.
The space $\mathcal D$ consists of all functions of the form
$$ f(x,s) = \sum_{k=0}^n s^k a_k (x),\quad (x,s)\in X\times\R,$$
where $n \in \mathbb{N}$ and $a_0(x),a_1(x),\dots, a_n(x) \in
\mathcal{D}(X)$. The nuclear space $\mathcal D$ is topologically embedded into the Hilbert space
\begin{equation}\label{fty7ed7i}
H_0 = L^2(\mathbb{R}^d \times\mathbb{R},dx\,d\sigma(s)).\end{equation}
Thus, we get a Gel'fand triple
$$ \mathcal{D} \subset H_0 \subset \mathcal{D}',$$
where $\mathcal {D}'$ is the dual space of $ \mathcal{D}$ with
respect to the zero space $H_0$.

For each $\varphi\in\mathcal D(X)$, we define $\varphi\otimes1, \varphi\otimes\operatorname{id}\in \mathcal D$ by
$$(\varphi\otimes1)(x,s)=\varphi(x),\quad  (\varphi\otimes\operatorname{id})(x,s)=\varphi(x)s,\quad (x,s)\in X\times\R.$$

Just as in Section~\ref{ghfdy7er5}, formula \eqref{huyf67e6},
 we construct the standard triple
$$\mathcal{F}_{\text {fin}}(\mathcal{D}) \subset\mathcal{F}(H_0) \subset \mathcal{F}^*_{\text
{fin}}(\mathcal{D}).$$

The following theorem was proved in \cite{Das}. It contains Theorems~\ref{buf7o6} and \ref{gfye6} as special cases.

\begin{theorem}
Set $\Phi=\mathcal D(X)$, $H=L^2(X,dx)$, $\Phi'=\mathcal D'(X)$, $\Psi=\mathcal{F}_{\mathrm {fin}}(\mathcal{D})$,
$\mathcal F=\mathcal F(H_0)$, $\Psi'=\mathcal{F}^*_{\mathrm {fin}}(\mathcal{D})$.
 For each $\varphi \in \Phi$, we define a linear operator
\begin{equation}\label{bvufr7o}
 A(\varphi )= a^+(\varphi \otimes 1) + a^0(\varphi\otimes\operatorname{id}) +a^-(\varphi \otimes 1) \end{equation}
acting on $\Psi$. Then conditions of Theorem~\ref{spectral} are satisfied. The corresponding measure $\mu$ is the L\'evy white noise measure with Kolmogorov measure $\sigma$:
 \begin{align}
 &\int_{\mathcal D'(X)}e^{i\langle\omega,\varphi\rangle }\,d\mu(\omega)
  =\exp\bigg[
 \int_{X\times\R}(e^{is\varphi(x)}-1-is\varphi(x))\frac1{s^2}\,dx\,d\sigma(s)
 \bigg],\quad \varphi\in\mathcal D(X).\label{hdfyrrd}\end{align}
 The function under the integral sign on the right hand side of formula \eqref{hdfyrrd} takes, for $s=0$, the limiting value
 $-(1/2)\varphi^2(x)$.
 \end{theorem}

Let us assume that $\sigma(\{0\})=0$ and the measure $\sigma$ additionally satisfies $\int_{\R\setminus\{0\}} \frac1s\,d\sigma(s)<\infty$. Instead of \eqref{bvufr7o}, we may consider the operators
$$ A(\varphi )= a^+(\varphi \otimes 1)+ a^0(\varphi\otimes\operatorname{id})+ a^-(\varphi \otimes 1) +\int_{X}
\varphi(x)\,dx\, \int_{\R\setminus\{0\}} \frac1s\,d\sigma(s).
$$
The corresponding measure $\mu$ has  the Fourier transform
$$\int_{\mathcal D'(X)}e^{i\langle\omega,\varphi\rangle }\,d\mu(\omega)=\exp\left[
 \int_{X\times(\R\setminus\{0\})}(e^{is\varphi(x)}-1)\frac1{s^2}\,dx\,d\sigma(s)
 \right],\quad \varphi\in\mathcal D(X).$$

 Let $\mathbb K(X)$ denote the set of all discrete signed Radon measures on $X$:
 $$\mathbb K(X)=\left\{\eta=\sum_i s_i\delta_{x_i}\mid s_i\in\R_+,\ x_i\in\R,\ \text{$\eta$ is a signed Radon measure on $X$}\right\}.$$
 Evidently, $\mathbb K(X)\subset\mathcal D'(X)$. Furthermore, it can be shown that $\mathbb K(X)\in\mathcal C(\mathcal D'(X))$. The measure $\mu$ is concentrated on $\mathbb K(X)$, i.e.,  $\mu(\mathbb K(X))=1$, cf.\ \cite{KLV}. Analogously, if the measure $\sigma$ is concentrated on $(0,\infty)$, then $\mu(\mathbb K_+(X))=1$. Here
 $$\mathbb K_+(X)=\{\eta\in\mathbb K(X)\mid \text{$\eta$ is a Radon measure on $X$}\}.$$

If the measure $\frac1{s^2}\,d\sigma(s)$ is finite, then $\mu$-a.s.\   the set of atoms of $\eta=\sum_{i}s_i\delta_{x_i}$ is locally finite in $X$, i.e., $\{x_i\}\in\Gamma(X)$. Thus, the measure $\mu$ is a marked Poisson process. In particular, we recover the Poisson measure $\mu$ when $\sigma=\delta_1$.

On the other hand, if the measure $\frac1{s^2}\,d\sigma(s)$ is infinite, then $\mu$-a.s.\  the set of atoms of $\eta=\sum_{i}s_i\delta_{x_i}$, i.e., the set $\{x_i\}$, is dense in $X$.

For example, let $\alpha,\beta>0$ and let
$$d\sigma(s)=\beta se^{-s/\alpha}\,ds.$$
Then $\mu$ is a gamma measure, cf.\ \cite{TsVY}. It can be easily calculated that its Laplace transform is given by
$$\int_{\mathbb K_+(X)}\exp[-\langle\eta,\varphi\rangle]\,d\mu(\eta)=\exp\left[-\beta\int_X\log(1+\alpha\varphi(x))\,dx\right],\quad \varphi \in\mathcal D(X),\ \varphi(x)>-1/\alpha.$$

\section{Further developments}\label{vgfvyf}

\subsection{Generalized stochastic processes with independent values}

 It is possible to consider the case where the measure $\sigma$ from Section~\ref{gyd6ei6} depends on $x$. That is for each $x \in \mathbb {R}^d$, $\sigma(x,ds)$ is a probability measure on
$(\mathbb{R},\mathcal{B}(\mathbb{R}))$. We also assume that for each $\Delta \in \mathcal{B}(\mathbb{R})$,
\begin{align}
\mathbb{R}^d \ni x\mapsto \sigma (x, \Delta)
\end{align}
is a measurable mapping.
(Note that, if $d=1$,  $\sigma(x,ds)$ is just a Markov kernel on $(\mathbb R,\mathcal B(\mathbb R))$.)
Hence,  we can define a $\sigma$-finite measure
$dx\, \sigma(x,ds)$ on $(\mathbb{R}^d \times \mathbb{R},\mathcal{B}(\mathbb{R}^d \times \mathbb{R}))$.
By analogy with estimate \eqref{equ 3.2}, one assumes that, for each locally compact $\Lambda\in\mathcal B(X)$, there exists a constant $C_\Lambda>0$ such that
\begin{equation}\label{gvfytf}
\int_\Lambda \int_\R   s^n \,\sigma(x,ds)\,dx
\le C_\Lambda^n\,n!,\quad n\in\mathbb N.\end{equation}
Then the family of operators $A(\varphi)$ of the form \eqref{bvufr7o} satisfies the conditions of Theorem~\ref{spectral} and the corresponding spectral measure $\mu$ has Fourier transform
$$\int_{\mathcal D'(X)}e^{i\langle\omega,\varphi\rangle }\,d\mu(\omega)=\exp\left[
 \int_{X} \int_\R(e^{is\varphi(x)}-1)\frac1{s^2}\,d\sigma(x,s)
 \,dx\right],\quad \varphi\in\mathcal D(X),$$
see \cite{Das} for further details. Thus, $\mu$ is a generalized stochastic process with independent values \cite{GV}.

For example, let $\alpha,\beta:X\to(0,\infty)$ be measurable functions. Furthermore, let us assume that the function $\alpha$ is locally bounded and $\alpha\beta\in L^1_{\mathrm{loc}}(X,dx)$. We define
$$ \sigma(x,ds)=\beta(x)s e^{-s/\alpha(x)}\,ds. $$
As easily seen, \eqref{gvfytf} holds. The  corresponding measure $\mu$ is concentrated on $\mathbb K_+(X)$ and has Laplace transform
\begin{multline*}
\int_{\mathbb K_+(X)}\exp[-\langle\eta,\varphi\rangle]\,d\mu(\eta)=\exp\left[-\int_X \beta(x)\log\big(1+\alpha(x)\varphi(x)\big)\,dx\right],\\
\varphi \in\mathcal D(X),\ \varphi(x)>-1/\alpha(x).\end{multline*}

\begin{remark}
 In some cases, for a given family of Hermitian operators $A(\varphi)$, $\varphi\in\Phi$, one cannot find a nuclear space $\Psi$ for which the conditions of Theorem~\ref{spectral} are satisfied. Instead, one has to use a linear topological space $\Psi$ which is topologically  embedded into a Hilbert space $\mathcal F$.  In such a case, one cannot derive existence of a corresponding spectral measure $\mu$ directly from the projection spectral theorem, so additional considerations are required.
Still, in papers \cite{L2,L3} determinantal (fermion) and permanental (boson) point processes were derived by using such an approach. This, in particular, showed that each determinantal point process is the spectral measure of a family of operators obtained by smearing out the particle density for a quasi-free representation of Canonical Anticommutation Relations.
\end{remark}

 \subsection{Free L\'evy white noise}

 A principal assumption of the projection spectral theorem is the commutation of the operators $A(\varphi)$. However, there are important situations when self-adjoint  operators $A(\varphi)$ which form a Jacobi field do not commute. So there is no corresponding spectral measure $\mu $. Still, in some cases, it is possible to study a noncommutative Jacobi field $(A(\varphi))_{\varphi\in\Phi}$ within the framework of noncommutative probability. In particular, there is a deep theory of L\'evy white noise in free probability. Let us very briefly mention some elements of this theory and refer  the reader to \cite{BL1,BL2} for a detailed discussion of it.

Let $H_0$ be a real separable Hilbert space and let $F(H_0)$ denote the full Fock space over $H_0$:
$$F(H_0)=\bigoplus_{n=0}^\infty H_0^{\otimes n}.$$
For each $\varphi\in H_0$ we define a creation operator $a^+(\varphi)$ and an annihilation operator $a^-(\varphi)$ as bounded linear operators on $F(H)$ satisfying
\begin{align*}
a^+(\varphi)f^{(n)}&=\varphi\otimes f^{(n)},\quad f^{(n)}\in H_0^{\otimes n},\\
a^-(\varphi)h_1\otimes h_2\otimes\dots\otimes h_n&=(\varphi,h_1)_{H_0}\,h_2\otimes\dots\otimes h_n, \quad h_1,h_2,\dots,h_n\in H_0.
\end{align*}
Assume that $H_0$ is an $L^2$-space and $\varphi\in H$ is a bounded function. Then, we define a neutral operator $a_0(\varphi)$
as a bounded linear operator on $F(H)$ satisfying
$$a^0(\varphi)h_1\otimes h_2\otimes\dots\otimes h_n=(\varphi h_1)\otimes h_2\otimes \dots\otimes h_n, \quad h_1,h_2,\dots,h_n\in H_0.$$

Let $\sigma$ be a finite, nonzero measure on $\R$ with compact support and let the Hilbert space $H_0$ be given by \eqref{fty7ed7i}.
For each $\varphi\in\mathcal D(X)$, we define a  self-adjoint bounded linear operator $A(\varphi)$ on $F(H_0)$ just as in  formula \eqref{bvufr7o}. Let $\mathbf A$ denote the real algebra generated by the operators $(A(\varphi))_{\varphi\in\mathcal D(X)}$. We define a free expectation on $\mathbf A$ by
$$\tau (\mathbf a)=(\mathbf a\Omega,\Omega)_{F(H_0)},\quad \mathbf a\in\mathbf A,$$
where $\Omega=(1,0,0,\dots)$ is the vacuum vector in $F(H_0)$.

Recall that a set partition $\pi$ of a set $Y$ is a collection of disjoint subsets of $Y$ whose union equals $Y$. Let $ {NC}(n)$ denote the collection of all non-crossing partitions of $\{1,\dots,n\}$, i.e., all set partitions $\pi=\{B_1,\dots,B_k\}$, $k\ge1$, of $\{1,\dots,n\}$ such that there do not exist $B_i,B_j\in\pi$, $B_i\ne B_j$, for which the following inequalities hold: $b_1<b_3<b_2<b_4$
for some $b_1,b_2\in B_i$ and $b_3,b_4\in B_j$.

For each $n\in\mathbb N$, we define a free cumulant $C^{(n)}$ as the $n$-linear mapping $C^{(n)}:\mathcal D(X)^n\to\mathbb R$
given recurrently by the following formula, which connects the free cumulants with moments:
\begin{equation}\label{jkgi8y}
\tau(A(\varphi_1)A(\varphi_2)\dotsm A(\varphi_n))=\sum_{\pi\in{NC}(n)}\prod_{B\in\pi}C(B,\varphi_1,\dots,\varphi_n),\end{equation}
where for each $B=\{b_1,\dots,b_k\}\subset \{1,2,\dots,n\}$, $b_1<b_2<\dots<b_k$,
$$ C(B,\varphi_1,\dots,\varphi_n):=C^{(k)}(\varphi_{b_1},\dots,\varphi_{b_k}). $$

It can be shown that, in our case $C^{(1)}\equiv0$ and
\begin{equation}\label{hjfdtyr} C^{(n)}(\varphi_1,\dots,\varphi_n)=\int_{\R}s^{n-2}\,d\sigma(s)\,\int_X \varphi_1(x)\dotsm \varphi_n(x)\,dx,
 \end{equation}
 where $\varphi_1,\dots, \varphi_n\in\mathcal D(X)$ and $n\ge2$.
Hence, by \eqref{jkgi8y} and \eqref{hjfdtyr}, the expectation $\tau$ on $\mathbf A$ is tracial, i.e., for any $\mathbf a,\mathbf b\in\mathbf A$, $\tau(\mathbf{ab})=\tau(\mathbf{ba})$.

Let $\varphi_1,\dots, \varphi_n\in\mathcal D(X)$ be such that $\varphi_i\varphi_j=0$ if $i\ne j$. Then, by \eqref{hjfdtyr}, for each $k\ge2$ and any indices $i_1,\dots,i_k\in\{1,\dots,n\}$ such that $i_l\ne i_m$ for some $l,m\in\{1,\dots,k\}$, $C^{(k)}(\varphi_{i_1},\dots,\varphi_{i_k})=0$.
This means that the operators $A(\varphi_1),\dots,A(\varphi_n)$ are freely independent with respect to the expectation $\tau$, see \cite{Speicher}.
Hence, if we introduce (at least informally) operators $A(x)$ ($x\in X$), so that $A(\varphi)=\int_X\varphi(x)A(x)\,dx$,
then we may think of $A(x)$ as a free L\'evy white noise.

Let us recall that, in classical probability, the cumulant transform of a probability measure can be expressed as the logarithm of the Laplace transform of the measure. In particular, for a probability measure $\mu$ on $\mathcal D'(X)$, the cumulant transform of $\mu$ is given by
$$C_\mu(\varphi)=\log\bigg[\int_{\mathcal D'(X)}e^{\langle\omega,\varphi\rangle}\,d\mu(\omega)\bigg],\quad \varphi\in\mathcal D(X).$$
The counterpart of the cumulant transform in free probability is the free cumulant transform. In our case, the free cumulant transform is defined by
\begin{equation}\label{kjlhiuy}
C(\varphi)=\sum_{n=1}^\infty C^{(n)}(\varphi,\dots,\varphi) \end{equation}
for $\varphi\in\mathcal D(X)$ such that the series on the right hand side of \eqref{kjlhiuy} converges.

It is shown in \cite{BL1} that there exists $\varepsilon>0$ such that, for each $\varphi\in\mathcal D(X)$ with $|\varphi|<\varepsilon$, we have
\begin{equation}\label{kjgyt}
C(\varphi)=\int_{X\times \R} \frac{\varphi^2(x)}{1-s\varphi(x)}\,dx\,d\sigma(s).\end{equation}
In particular, for the free Gaussian white noise, we have $\sigma=\delta_0$ and so
$$ C(\varphi)=\int_{X} \varphi^2(x)\,dx.$$
Note, for comparison, that for the Gaussian white noise measure $\mu$ from Theorem~\ref{buf7o6}, we have
$$C_\mu(\varphi)=\frac12\int_{X} \varphi^2(x)\,dx.$$

For the free Poisson white noise, we have $\sigma=\delta_1$ and so
$$C(\varphi)=\int_{X} \frac{\varphi^2(x)}{1-\varphi(x)}\,dx=\int_X\sum_{n=2}^\infty \varphi^{n}(x)\,dx.$$
Again, for comparison, for the centered Poisson measure $\mu$, we have
$$C_\mu(\varphi)=\int_X\big(e^{\varphi(x)}-1-\varphi(x)\big)\,dx=\int_X\sum_{n=2}^\infty \frac{\varphi^n(x)}{n!}\,dx.$$

Finally, for a general free L\'evy white noise, we get, by \eqref{kjgyt},
$$C(\varphi)=\sum_{n=2}^\infty\int_\R s^{n-2}\,d\sigma(s)\int_X\varphi^n(x)\,dx,$$
while for a classical L\'evy white noise $\mu$, we get by \eqref{hdfyrrd},
$$C_\mu(\varphi)=\sum_{n=2}^\infty\frac1{n!}\int_\R s^{n-2}\,d\sigma(s)\int_X\varphi^n(x)\,dx.$$

\begin{center}
{\bf Acknowledgements}\end{center}
The author acknowledges the financial support of the SFB 701 ``Spectral
structures and topological methods in mathematics'', Bielefeld University and of the Polish
National Science Center, grant no.\ Dec-2012/05/B/ST1/00626.


\begin{thebibliography}{99}

\bibitem{ber_1968} Berezanskii, Ju.M.: Expansions in eigenfunctions of selfadjoint operators. American Mathematical Society, Providence, R.I. 1968.

\bibitem{ber_1984} Berezanskii, Yu.M.: The projection spectral theorem. Russian Math. Surveys 39 (1984), 1--62.

\bibitem{ber_1985} Berezanskii, Yu.M.: On the projection spectral theorem. Ukrainian Math. J. 37 (1985), 124--130.

\bibitem{ber_1978}
Berezanskii, Yu.M.:  Selfadjoint operators in spaces of functions of infinitely many variables. Naukova Dumka, Kiev, 1978 (in Russian), English translation: American Mathematical Society, Providence, RI, 1986.


\bibitem{ber_1993} Berezansky, Yu.M.: Spectral approach to white noise analysis. Dynamics of complex and irregular systems (Bielefeld, 1991), 131--140, Bielefeld Encount. Math. Phys., VIII, World Sci. Publ., River Edge, NJ, 1993.

\bibitem{ber_2000} Berezansky, Y.M.:
Poisson measure as the spectral measure of Jacobi field. Infin. Dimens. Anal. Quantum Probab. Relat. Top. 3 (2000),  121--139.


\bibitem{BK} Berezansky, Y.M., Kondratiev, Y.G.: Spectral methods in infinite-dimensional analysis. Vol. 1, 2.  Kluwer Academic Publishers, Dordrecht, 1995.


\bibitem{BKosh}  Berezanskii, Ju.M.,  Koshmanenko, V.D. Axiomatic field theory in terms of operator Jacobi matrices. (Russian) Teoret. Mat. Fiz. 8 (1971), 175--191.

\bibitem{BLL}  Berezanskii, Yu.M., Livinskii, V.O., Lytvynov, E.W.: A spectral approach to white  noise analysis. (Russian) Ukrain. Mat. Zh. 46 (1994), 177--197; translation in Ukrainian Math. J. 46 (1994), 183--203 (1995).

\bibitem{BL1}     Bo\.zejko, M., Lytvynov, E.:  Meixner class of non-commutative generalized stochastic processes with freely independent values. I. A characterization. Comm. Math. Phys. 292 (2009),  99--129.

\bibitem{BL2}     Bo\.zejko, M., Lytvynov, E.:  Meixner class of non-commutative generalized stochastic processes with freely independent values II. The generating function. Comm. Math. Phys. 302 (2011), 425--451.


\bibitem{Das} Das, S.: Orthogonal decompositions for generalized  stochastic processes with independent values, PhD dissertation, Swansea University, 2012.

\bibitem{GV}  Gel'fand, I. M., Vilenkin, N.Y.: Generalized functions. Vol. 4: Applications of harmonic analysis. Academic Press, New York-London, 1964.

\bibitem{HKPS} Hida, T.,  Kuo, H.-H., Potthoff, J., Streit, L.: White noise: An
  infinite dimensional calculus. Dordrecht, Boston,
  London: Kluwer Acad.\ Publ.,  1993.

\bibitem{IK} Ito, Y., Kubo, I.:
Calculus on Gaussian and Poisson white noises.
Nagoya Math. J. 111, 41--84 (1988).

\bibitem{Kallenberg} Kallenberg, O.:  Random measures. Akademie-Verlag, Berlin; Academic Press, London-New York, 1976.

\bibitem{KLV} Kondratiev, Y., Lytvynov, E, Vershik, A.:  Laplace operators on the cone of Radon measures, in preparation.

\bibitem{L1}  Lytvynov, E.W.: Multiple Wiener integrals and non-Gaussian white noises: a Jacobi field approach. Methods Funct. Anal. Topology 1 (1995), no. 1, 61--85.

\bibitem{L2} Lytvynov, E.:  Fermion and boson random point processes as particle distributions of infinite free Fermi and Bose gases of finite density. Rev. Math. Phys. 14 (2002),  1073--1098.

\bibitem{L3} Lytvynov, E.:  Determinantal point processes with J-Hermitian correlation kernels. Ann. Probab. 41 (2013), 2513--2543.



    \bibitem{Speicher}  Speicher, R.: Free probability theory and non-crossing partitions.
 S\'em. Lothar.\ Combin.\  39, Art. B39c, 38 pp. (1997)  (electronic)

 \bibitem{TsVY} Tsilevich, N, Vershik, A,  Yor, M.: An infinite-dimensional analogue of the Lebesgue measure and distinguished properties of the gamma process. J. Funct. Anal. 185 (2001), 274--296.

\end{thebibliography}
\end{document}